%% file: jssc-module.tex
\documentclass{jssc}

\usepackage{graphicx}         
\usepackage{amsmath}
\allowdisplaybreaks
\usepackage{mathtools}
\usepackage{xcolor}
\usepackage{amssymb} 
\usepackage{setspace}
\usepackage{comment}
\usepackage{subfig}
\usepackage{float}
\usepackage[labelsep=newline,singlelinecheck=false,skip=0pt]{caption}
 \usepackage{xfrac,booktabs}
\usepackage{enumitem}
\usepackage[title]{appendix}

\DeclareMathOperator{\rank}{rank}
\DeclareMathOperator{\nullspace}{null}
\DeclareMathOperator{\vectorize}{vec}
\DeclareMathOperator{\halfvectorize}{vech}
\DeclareMathOperator{\matrixize}{mat}

\def\textsubscript#1%
{$_{\text{#1}}$}

\def\cdd{\mbox{\boldmath$\cdot$}~}

\newcommand{\rulex}{\hfill\rule{1mm}{3mm}}

\input{amssym.def}
\include{graphix}

\input jsscN.tex
\begin{document}

\title{Inverse Continuous-Time Linear Quadratic Regulator: From Control Cost Matrix to Entire Cost Reconstruction}
{\uppercase{CAO} Yuexin \cdd \uppercase{LI}
Yibei \cdd \uppercase{Zou} Zhuo \cdd \uppercase{HU} Xiaoming} 
{\uppercase{CAO} Yuexin \cdd \uppercase{HU} Xiaoming  \\
Department of Mathematics, KTH Royal Institute of Technology, Stockholm $10044$, Sweden.  \\ Email: yuexin@kth.se; hu@math.kth.se \\  
\uppercase{LI} Yibei    \\
State Key Laboratory of Mathematical Sciences, Academy of Mathematics and Systems Science, Chinese Academy of Sciences, Beijing $100190$, China.  Email: yibeili@amss.ac.cn \\
\uppercase{ZOU} Zhuo    \\ School of Information Science and Technology, Fudan University, Shanghai $200433$, China \\
  Email: zhuo@fudan.edu.cn
   } 
{
{$^\diamond$}}


\dshm{Year}{Volume}{\uppercase{Inverse Optimal Control}}{\uppercase{CAO Yuexin} \cdd \uppercase{LI
Yibei} \cdd \uppercase{Zou Zhuo} \cdd \uppercase{HU Xiaoming}} 



\Abstract{This paper studies the inverse optimal control problem for continuous-time linear quadratic regulators over finite-time horizon, aiming to reconstruct the control, state, and terminal cost matrices in the objective function from observed optimal inputs. Previous studies have mainly explored the recovery of state cost matrices under the assumptions that the system is controllable and the control cost matrix is given. Motivated by various applications in which the control cost matrix is unknown and needs to be identified, we present two reconstruction methods. The first exploits the full trajectory of the feedback matrix and establishes the necessary and sufficient condition for unique recovery. To further reduce the computational complexity, the second method utilizes the feedback matrix at some time points, where sufficient conditions for uniqueness are provided. Moreover, we study the recovery of the state and terminal cost matrices in a more general manner. Unlike prior works that assume controllability, we analyse its impact on well-posedness, and derive analytical expressions for unknown matrices for both controllable and uncontrollable cases. Finally, we characterize the structural relation between the inverse problems with the control cost matrix either to be reconstructed or given as a prior.}      



\section{Introduction}
Inverse optimal control (IOC) seeks to reconstruct the objective function based on the knowledge of the system dynamics and the observed optimal actions taken by the systems. Initially introduced by Kalman \citep{kalman1963linear}, IOC has received considerable attention and has been widely applied across various fields, including control engineering \citep{maillotpilotsfly2013,yokoyama2018inference}, robotics \citep{chittaro2013inverse,finn2016guided}, and economics \citep{saez2016datapriceelectricity}. The study of IOC problems is of great importance, as it enables the reconstruction of optimal control models and provides deeper insights into the underlying principles governing complex physical systems.

Depending on the application context, IOC problems can be formulated into different forms. One well-studied direction is to represent the objective function as a linear combination of some known basic functions, with the corresponding weights remaining to be identified \citep{Keshavarz2011imputconvex,hatz2012estimating,johnson2013inversebasicfunction,molloy2018finitenonlinear}. As a fundamental criterion in optimal control, IOC for linear quadratic regulator (LQR) has also received considerable attention. The objective function of the standard LQR is defined by three key cost matrices: the control cost matrix $R$, the state cost matrix $Q$ and the terminal cost matrix $F$ \citep{zabczyk2020mathematicalcontroltheory}. The recovery of cost matrices has been explored extensively for discrete-time systems over both infinite-time horizons \citep{priess2014solutions,molloy2020onlineIOC} and finite-time horizons \citep{zhang2019inverse,yu2021systemidentification,zhang2023convexapproach,zhang2024statistically}. For continuous-time systems over infinite-time horizons, related studies can be found in \citep{sugimoto1987newsolutionIOC,molloy2022inversebook}.

In contrast, significantly less attention has been devoted to IOC for continuous-time systems over finite-time horizons. Unlike the standard LQR, IOC problems studied in \citep{nori2004linearoptimalcrossterm,jean2018inversecrossterm} include a terminal state constraint and a state-control cross-term in the objective function, with solutions characterized by algebraic Riccati equations. Due to the involvement of differential Riccati equations, IOC problems with the standard LQR formulation are particularly challenging. In \citep{jameson1973inverse}, necessary and sufficient conditions are established for reconstructing the objective function based on the optimal feedback matrix. However, the resulting cost matrices may lack desirable properties: $R$ may not be constant, and $Q$ may be neither constant nor positive semi-definite. To address the limitations, \citep{li2020continuous} establishes the necessary and sufficient condition for the existence of cost matrices by using the framework of Linear Matrix Inequality (LMI). Furthermore, \citep{li2021identifiability} extends the analysis to scenarios with partial state observations under the assumption of system invertibility. Most of existing works focus on recovering constant and positive semi-definite $Q$ while assuming that $R$ is known and the underlying system is controllable.

In this paper, we study IOC problems for the standard continuous-time LQR over finite-time horizons with the optimal feedback matrix given. Previous work have primarily focused on the reconstruction of $Q$ and $F$ under the assumptions that $R$ is given and the underlying system is controllable. Motivated by the various applications in which $R$ is unknown, we first focus on its reconstruction. Then we consider IOC problems where $R$ is given in a more general setup. The contribution of this paper is threefold. Firstly, we propose two methods to recover $R$. The first method leverages the full trajectory of the optimal feedback matrix. We provide the necessary and sufficient condition for the unique recovery of $R$, and when this condition is not satisfied, we derive its solution space. Although the first method offers a detailed characterization of the solution space of $R$, it may be data-demanding and computationally expensive in certain cases. As a complement, we present the second method that recovers $R$ using the optimal feedback matrix at some time points. This method reduces computational burden and can be utilized to the cases where only $F$ is available or where data is collected in real time. Secondly, we consider the recovery of $Q$ and $F$ with the given $R$ in a more comprehensive manner. In contrast to prior studies that impose controllability as a standing assumption, we analyse its role in determining the well-posedness of the problem. When $F$ is further specified, we characterize the solution space of $Q$ in both of the controllable and uncontrollable cases, and vice versa. When $Q$ and $F$ are unknown, we derive the sufficient conditions for well-posedness and provide their analytical expressions. Thirdly, by introducing an auxiliary optimal control problem, we establish the structural connection between the inverse problem of recovering $R$ and the one in which $R$ is prescribed as a prior.


This paper is organized as follows: In Section \ref{sec:problem_formulation}, two IOC problems for continuous-time LQR over finite-time horizons are formulated. Section \ref{sec:reconstruction_R} studies the IOC problem where $R$ is unknown and presents two methods for its recovery. Section \ref{sec:Rknown} considers the case where $R$ is given, and focuses on recovering $Q$ and $F$. In Section \ref{sec:equivalence_p1_p2}, we analyze the structural connection between these two IOC problems. Section \ref{sec:conclusions_futurework} concludes the paper and outlines directions for future work.

\textit{Notations}: For any positive integer $n$ and $m$, let $\mathbb{R}$, $\mathbb{R}^{n}$ and $\mathbb{R}^{n \times m}$ denote the set of real numbers, $n \times 1$ real vectors and $n \times m$ real matrices, respectively. Let $\mathbb{S}^{n}_{+}$ and $\mathbb{S}^{n}_{++}$ denote the set of $n \times n$ positive semi-definite matrices and positive definite matrices, respectively. We denote $[n]$ as the set of natural numbers $\{1,2,\ldots,n\}$. Let $I$ denote the identity matrix whose dimensions are inferred from the context. Let $\otimes$ denote the Kronecker product.  For a matrix $A_{1}$, denote $\rank(A_{1})$, $\nullspace(A_{1})$ and $A_{1}^{\dagger}$ as the rank of $A_{1}$, the null space of $A_{1}$ and the Moore-Penrose inverse of $A_{1}$, respectively. For a square matrix $A_{2}$, let $\lambda_{\min}(A_{2})$ denote the minimal eigenvalue of $A_{2}$. Let $A_{3}(t)$ be a matrix whose elements are functions of the variable $t$, for any integer $i$, its $i$-th derivative is denoted by $A_{3}(t)^{(i)}$. We use $\vectorize(\cdot)$, $\halfvectorize(\cdot)$ and $\matrixize(\cdot)$ to denote vectorization, half-vectorization and matricization, respectively.


\section{Problem Formulation}\label{sec:problem_formulation}

The forward linear quadratic optimal control problems over finite-time horizons are formulated as follows:
\begin{equation}\label{equ:classic_optimal_problem}
\begin{split}
\min_{u}  \quad &  x^{\rm T}(t_{f})Fx(t_{f}) + \int_{t_{0}}^{t_{f}} \! (x^{\rm T}(t) Qx(t) +u^{\rm T}(t) Ru(t) )  \, \mathrm{d}t  \\
\text{s.t.} \quad & \dot{x}(t) = Ax(t) + Bu(t), \, x(t_{0}) = x_{0}, 
\end{split}
\end{equation}
where $x(t) \in \mathbb{R}^{n}$, $u(t) \in \mathbb{R}^{m}$, $Q, F \in \mathbb{S}^{n}_{+}$ and $R \in \mathbb{S}^{m}_{++}$. The underlying linear system, denoted by $(A, B)$, is assumed to be given and $m \leq n$. We denote $\rank(B)$ as $r_{b}$.

For the forward optimal control problem (\ref{equ:classic_optimal_problem}), there exists a unique optimal feedback controller:
\begin{equation}\label{equ:classic_forward_feedback_control}
u(t) = K(t)x(t) = -R^{-1}B^{\rm T}P(t)x(t),
\end{equation}
where the matrix $P(t)$ is the solution of the following differential Riccati equation (DRE) over the time interval $t \in [t_{0}, t_{f}]$:
\begin{equation*}
\begin{split}
\dot{P}(t) & = - A^{\rm T}P(t) - P(t)A + P(t)BR^{-1}B^{\rm T}P(t) - Q, \\ P(t_{f}) & = F.
\end{split}
\end{equation*}
It is well-known that if $Q, F \in \mathbb{S}^{n}_{+}$ and $R \in \mathbb{S}^{m}_{++}$, (DRE) admits a unique solution $P(t) \in \mathbb{S}^{n}_{+}$ over $[t_{0}, t_{f}]$.

In this paper, we assume that the optimal feedback matrix $K(t)$ is given. The goal of IOC problems is to recover $R$, $Q$ and $F$ in the objective function such that it could generate the observed optimal feedback controller. By (\ref{equ:classic_forward_feedback_control}), it holds that $B^{\rm T}P(t) = -RK(t)$ for $t \in [t_{0}, t_{f}]$. Based on the analysis of matrix equations, the necessary and sufficient conditions to reconstruct $R$ and $P(t)$ are provided in \citep{jameson1973inverse}, which are summarized in Lemma \ref{lemma:necessary_sufficient_conditions_1973}.

\begin{lemma}\label{lemma:necessary_sufficient_conditions_1973}
Given $K(t)$, the equation $B^{\rm T}P(t) = -RK(t)$ has solutions $R \in \mathbb{S}^{m}_{++}$ and $P(t) \in \mathbb{S}^{n}_{+}$ if and only if the following conditions hold for all $t \in [t_{0}, t_{f}]$:
\begin{enumerate}[label={(A\arabic*})]
\item \label{assumption1} $K(t)B$ has $m$ linearly independent real eigenvectors;
\item \label{assumption2} all the eigenvalues of $K(t)B$ are non-positive;  
\item \label{assumption3} $\rank(K(t)B) = \rank(K(t))$.
\end{enumerate}
To ensure $P(t) \in \mathbb{S}^{n}_{++}$, the third condition \ref{assumption3} is further strengthened to $\rank(K(t)B) = \rank(K(t)) = \rank(B)$.
\end{lemma}

Since it is required that $Q, F \in \mathbb{S}^{n}_{+}$ and $R \in \mathbb{S}^{m}_{++}$, we assume that \ref{assumption1}, \ref{assumption2} and \ref{assumption3} in Lemma \ref{lemma:necessary_sufficient_conditions_1973} are satisfied throughout this paper. When $K(t)$ is given without additional information, IOC problems are trivially ill-posed since proportional objective functions could generate the same optimal controller, namely $(aQ, aF, aR)$ will generate the same $K(t)$ for any $a > 0$. To eliminate this ambiguity and simplify IOC problems, most of the studies assume that $R = I$. In this paper, since $R$ is unknown in various applications and needs to be identified, we first focus on its recovery, which is stated in Problem \ref{problem:one_reconstruct_R}.

\begin{problem}\label{problem:one_reconstruct_R}
Given the optimal feedback matrix $K(t)$, recover the control cost matrix $R$ that governs the optimal behavior of the system when 
\begin{enumerate}[label={(1\alph*)}] 
    \item \label{problem:1a} the matrices $Q$ and $F$ are both known,
    \item \label{problem:1b} only the matrix $F$ is known.
\end{enumerate}
\end{problem}

Furthermore, IOC problems with the given $R$ are investigated in a more comprehensive manner. The common controllability assumption is relaxed. Since it holds that $RK(t)B = -B^{\rm T}P(t)B$ by (\ref{equ:classic_forward_feedback_control}), the symmetry of $P(t)$ implies that $R \in \mathcal{R}_{s} $, where $\mathcal{R}_{s}$ is defined as
\begin{equation*}
\mathcal{R}_{s} = \{R \in \mathbb{S}^{m}_{++} \mid RK(t)B = B^{\rm T}K(t)^{\rm T}R, \, \forall t \in [t_{0}, t_{f}] \}.
\end{equation*}
Then the reconstruction of $Q$ and/or $F$ are formulated in Problem \ref{problem:two_R_given}.

\begin{problem}\label{problem:two_R_given}
Given the optimal feedback matrix $K(t)$ and $R \in \mathcal{R}_{s}$, recover the cost matrices that govern the optimal behavior of the system under the following scenarios:
\begin{enumerate}[label={(2\alph*)}]
    \item \label{problem:2a} if $Q$ (or $F$) is further given, determine the corresponding $F$ (or $Q$);
    \item \label{problem:2b} If no additional prior knowledge is available, recover the matrices $Q$ and $F$.
\end{enumerate}
\end{problem}

In the following, Problem \ref{problem:one_reconstruct_R} and Problem \ref{problem:two_R_given} will be investigated in Section \ref{sec:reconstruction_R} and Section \ref{sec:Rknown}, respectively. In Section \ref{sec:equivalence_p1_p2}, the structural connection between these two problems will be further discussed.


\section{Reconstruction of $R$}\label{sec:reconstruction_R}

In this section, we investigate the problems of recovering $R$. We begin by analysing Problem \ref{problem:1a} and present two methods. The first method evaluates $R$ based on the full trajectory of $K(t)$ over $[t_{0}, t_{f}]$. The necessary and sufficient condition for the well-posedness is provided. When this condition is not satisfied, we characterize its solution space. As an alternative, we introduce a second method that recovers $R$ using $K(t)$ at some time points. This approach can be naturally extended to the cases where only $F$ is available, i.e., Problem \ref{problem:1b} or to real-time data collection scenarios, in which the first method is not applicable. Based on the second method, we also provide sufficient conditions for the unique recovery of $R$ and derive its analytical expression.

When $Q, F \in \mathbb{S}^{n}_{+}$ are given, the unique solution $P(t)$ of (DRE) can be explicitly computed, even though $R$ is unknown. By $K(t) = -R^{-1}B^{\rm T}P(t)$, (DRE) can be reformulated as the following Sylvester differential equation:
\begin{equation}\label{equ:matrix_DRE_without_R}
\begin{split}
\dot{P}(t) & = -A^{\rm T}P(t) - P(t)(A + BK(t)) - Q, \\ 
P(t_{f}) & = F.
\end{split}
\end{equation}
Then the unique $P(t)$ is obtained by solving (\ref{equ:matrix_DRE_without_R}), whose solution is given in \citep{abou2012matrixriccatiequation}. Since $P(t) \in \mathbb{S}^{n}_{+}$, then $\dot{P}(t)$ must be symmetric, i.e., $A^{\rm T}P + PA + Q + PBK = PA + A^{\rm T}P + Q + (PBK)^{\rm T}$,  
which implies that $P(t)BK(t)$ is symmetric for $t \in [t_{0}, t_{f}]$.

Consider Problem \ref{problem:1a}. We introduce the first method to recover $R$ by exploiting the full trajectory of $K(t)$. In the following theorem, we provide the necessary and sufficient condition for its well-posedness, and show that the full column rank of the input matrix $B$ is a necessity. When the condition does not hold, we derive the solution space of $R$.

\begin{theorem}\label{theorem:GivenFQ_R_wellposedness}
Given $K(t)$ and $Q, F \in \mathbb{S}^{n}_{+}$, the matrix $R$ satisfying $B^{\rm T}P(t) = -RK(t)$ is unique if and only if all row vectors of $K(t)$ are linearly independent over $[t_{0}, t_{f}]$, and the unique solution is given as $R = L_{1}^{-1}L_{2}$, where
\begin{equation*}\label{equ:integral_unique_R}
L_{1} = \int_{t_{0}}^{t_{f}} \! K(s)K(s)^{\rm T} \, \mathrm{d}s, \quad L_{2} = \int_{t_{0}}^{t_{f}} \! -K(s)P(s)B \, \mathrm{d}s.
\end{equation*}
A necessary condition for well-posedness is that the matrix $B$ has full column rank. 

If row vectors of $K(t)$ are linearly dependent, all the symmetric solutions are given as $R = \bar{R} +  V_{1}ZV_{1}^{\rm T}$, where $Z = Z^{\rm T}$, $V_{1}$ is the orthonormal basis for $\nullspace(L_{1})$ and \begin{equation}\label{equ:particular_solution_full_R}
\bar{R} = L_{1}^{\dagger}L_{2} + L_{2}^{\rm T}L_{1}^{\dagger} - L_{1}^{\dagger}L_{1}L_{2}^{\rm T}L_{1}^{\dagger}.
\end{equation}   
Then the solution space of $R$ is given as $\mathcal{U}_{R} = \{R \mid R = \bar{R} +  V_{1}ZV_{1}^{\rm T}\} \cap \mathbb{S}^{m}_{++}$.

Furthermore, if $\bar{R}$ in (\ref{equ:particular_solution_full_R}) is positive definite, then all the positive definite matrix $R$ can be explicitly characterized as $R = \bar{R} + V_{1}ZV_{1}^{\rm T}$, with $Z = Z^{\rm T}$ satisfying 
\begin{equation}\label{equ:condition_Z_R_pd}
Z \succ (V_{1}^{\rm T}\bar{R}V_{2})(V_{2}^{\rm T}\bar{R}V_{2})^{-1}(V_{2}^{\rm T}\bar{R}V_{1}) - V_{1}^{\rm T}\bar{R}V_{1},
\end{equation} 
where $V_{2}$ is the orthonormal basis for the range  of $L_{1}$.

\end{theorem}

\proof Firstly, by $B^{\rm T}P(t) = -RK(t)$, the following equation holds for all $t \in [t_{0}, t_{f}]$:
\begin{equation*}
\int_{t_{0}}^{t_{f}} \! K(s)K(s)^{\rm T} \, \mathrm{d}s R = \int_{t_{0}}^{t_{f}} \! -K(s)P(s)B \, \mathrm{d}s,
\end{equation*}
then the matrix $R$ has a unique solution if and only if $K(t)$ has linearly independent rows over $[t_{0}, t_{f}]$, and the unique solution is given as $R = L_{1}^{-1} L_{2}$.

Secondly, we establish that the condition for $B$ to have full column rank is necessary for the well-posedness. In the case where $B$ does not have full column rank, there exists a constant vector $w$ such that $Bw = 0$. Denote $R^{\star}$ as the true value of the control cost matrix. Then there exists a constant vector $v = R^{\star}w$ such that $B(R^{\star})^{-1}v = 0$, which indicates that $P(t)B(R^{\star})^{-1}v = 0$. Then the vector $v$ satisfies $K(t)^{\rm T} v = 0$, which implies that $K(t)$ does not have full row rank over $t \in [t_{0}, t_{f}]$.  

Next, we discuss the case where $K(t)$ has linearly dependent row vectors. Since there exists a solution $R^{\star}$, then $L_{1}R=L_{2}$ must be consistent, i.e., $L_{1}L_{1}^{\dagger}L_{2} = L_{2}$. When $L_{1}$ is singular, it holds for $R = L_{1}^{\dagger}L_{2}$, but it is not necessarily symmetric. Our constructed $\bar{R}$ in (\ref{equ:particular_solution_full_R}) is symmetric since $L_{1}L_{2}^{\rm T} = L_{1}RL_{1}$ and $L_{1} = L_{1}^{\rm T}$. It also holds that
\begin{equation*}
\begin{split}
L_{1}\bar{R} = L_{1}L_{1}^{\dagger}L_{2} + L_{1}L_{2}^{\rm T}L_{1}^{\dagger} - L_{1}L_{1}^{\dagger}L_{1}L_{2}^{\rm T}L_{1}^{\dagger} = L_{2},
\end{split}
\end{equation*}
which shows that $\bar{R}$ is symmetric. Then all the symmetric solutions are given as $\bar{R} + V_{1}ZV_{1}^{\rm T}$. The solution space $\mathcal{U}_{R}$ follows to ensure the positive definiteness. 

Finally, if $\bar{R}$ in (\ref{equ:particular_solution_full_R}) is positive definite, an explicit condition on $Z$ to ensure the positive definiteness of $R$ can be derived. Since $L_{1}$ is symmetric, the orthogonal complement of $\nullspace(L_{1})$ is the range space of $L_{1}$, then $U_{1} \coloneqq \begin{bmatrix} V_{1} & V_{2}  \end{bmatrix}$ is an orthonormal basis for $\mathbb{R}^{m}$. One can rewrite $R = \bar{R} + V_{1}ZV_{1}^{\rm T}$ in the basis of $U_{1}$ as:
\begin{equation*}
R = U_{1}\begin{bmatrix} V_{1}^{\rm T}\bar{R}V_{1} + Z & V_{1}^{\rm T}\bar{R}V_{2} \\ V_{2}^{\rm T}\bar{R}V_{1} & V_{2}^{\rm T}\bar{R}V_{2}\end{bmatrix} U_{1}^{\rm T} \coloneqq U_{1} U_{2} U_{1}^{\rm T}.
\end{equation*}
Since $U_{1}$ is non-singular, then $R \in \mathbb{S}^{m}_{++}$ if and only if $U_{2} \in \mathbb{S}^{m}_{++}$. The matrix $V_{2}^{\rm T}\bar{R}V_{2} \succ 0$ since $V_{2}$ has full column rank. Then by Schur complement test, the matrix $R \in \mathbb{S}^{m}_{++}$ if and only if $Z$ satisfies (\ref{equ:condition_Z_R_pd}). In fact, when $\bar{R}$ is positive definite, let $V_{1} \in \mathbb{R}^{m \times k}$, where $0< k < m$. The right-hand side of (\ref{equ:condition_Z_R_pd}) is a $k \times k$ matrix. Then the whole solution space can be characterized more efficiently, especially in cases where $k \ll m$.  \rulex 

\begin{remark}
Our first method exploits the entire trajectory of $K(t)$. If $L_{1}$ is non-singular, the matrix $R$ can be uniquely recovered by $R = L_{1}^{-1}L_{2}$. If $L_{1}$ is singular and a particular solution $\bar{R} \succ 0$ is available, the solution space of $R$ can be readily derived with an explicit condition on $Z$. In the absence of a positive definite particular solution, one can characterize the solution space as  $\mathcal{U}_{R}$.
\end{remark}


While the first method provides a detailed reconstruction of the solution space of $R$, it may be data-demanding and computationally intensive in certain cases. It also requires solving $P(t)$ for all $t \in [t_{0}, t_{f}]$, making it inapplicable to Problem \ref{problem:1b}. To address the limitations, we present a second method that uses $K(t)$ at some time points to recover $R$. This approach simplifies data requirements and reduces computational burden. Moreover, it can also be applied to the case where only $F$ is given or where data is received in real time. The second method will serve as a useful complement for characterizing the solution space of $R$.

Before we derive the solution space of $R$ by the second method, some properties of $BK(t)$ and useful facts about the orthogonal complement are firstly established in Proposition \ref{proposition:properties_BK} and Lemma \ref{lemma:space_direct_sum}, respectively.

\begin{proposition}\label{proposition:properties_BK}
For $\tilde{t} \in [t_{0}, t_{f}]$, let $r_{\tilde{t}}:=\rank(K(\tilde{t})B)$. If $1 \leq r_{\tilde{t}} \leq m$, denote negative eigenvalues of $K(\tilde{t})B$ as $\lambda_{\tilde{t}i}$ and its corresponding eigenvectors as $v_{\tilde{t}i}$, i.e., $K(\tilde{t})Bv_{\tilde{t}i} = \lambda_{\tilde{t}i}v_{\tilde{t}i}$, where $i \in [r_{\tilde{t}}]$. Then the following statements are valid:
\begin{enumerate}
\item $BK(\tilde{t})Bv_{\tilde{t}i} = \lambda_{\tilde{t}i}Bv_{\tilde{t}i}$, where $i \in [r_{\tilde{t}}]$;
\item $BK(\tilde{t})w_{\tilde{t}j} = 0$, where $K(\tilde{t})w_{\tilde{t}j} = 0$ and $j \in [n-r_{\tilde{t}}]$;
\item $\rank(BK(\tilde{t})) = \rank(K(\tilde{t})B) = r_{\tilde{t}}$.
\end{enumerate}
\end{proposition}

\proof See Appendix \ref{appendix:proof_BK}.  \rulex


\begin{lemma}\label{lemma:space_direct_sum}
Let $A \in \mathbb{R}^{m \times n}$. Denote $\mathbb{U}$ as its range space. Its orthogonal complement is $\mathbb{U}^{\perp} = \nullspace(A^{\rm T})$. For every $x \in \mathbb{R}^{m}$, there exist unique vectors $u_{1} \in \mathbb{U}$ and $u_{2} \in \mathbb{U}^{\perp}$ such that $x = u_{1} + u_{2}$.
\end{lemma}

Inspired by \citep{jameson1973inverse}, we present the analytical expression of $R$ when $P(t) \in \mathbb{S}^{n}_{++}$ for $t \in [t_{0}, t_{f}]$ in the next theorem.

\begin{theorem}\label{theorem:givenP_findR}
For a time point $\tilde{t} \in [t_{0}, t_{f}]$ such that $P(\tilde{t})$ is a positive definite matrix, all the positive definite $R(\tilde{t})$ satisfying $B^{\rm T}P(\tilde{t}) = -R(\tilde{t})K(\tilde{t})$ are represented by
\begin{equation}\label{equ:matrixR_Ppositivedefinite}
R_{1}(\tilde{t}) = -B^{\rm T}P(\tilde{t})(P(\tilde{t})BK(\tilde{t}))^{\dagger}P(\tilde{t})B + WW^{\rm T}, 
\end{equation}
where $W \in \mathbb{R}^{m \times (m-r_{b})}$ has linearly independent columns and $K(\tilde{t})^{\rm T}W = 0$. 
\end{theorem}

\proof In the following proof, our analysis is conducted on the time point $\tilde{t}$, allowing us to omit $'(\tilde{t})'$ for clarity and brevity. Since $P$ is positive definite, Lemma \ref{lemma:necessary_sufficient_conditions_1973} implies that $\rank(KB) = r_{b}$. By Proposition \ref{proposition:properties_BK}, one can obtain that $\rank(BK) = r_{b}$ and $BK$ is negative semi-definite. Since $PBK$ is symmetric, then $PBK$ is also negative semi-definite




We first show that $R_{1}$ in (\ref{equ:matrixR_Ppositivedefinite}) satisfies $B^{\rm T}P = -RK$ when $P \in \mathbb{S}^{n}_{++}$. It holds that $K^{\rm T}B^{\rm T}[P(PBK)^{\dagger}PBK-P] = 0$. Since $\rank(BK) = r_{b}$, by rank-nullity theorem, the null space of $K^{\rm T}B^{\rm T}$ has the same dimension as that of $B^{\rm T}$. Then $\nullspace(K^{\rm T}B^{\rm T}) = \nullspace(B^{\rm T})$ follows from the fact that $\nullspace(B^{\rm T})$ is a subspace of $\nullspace(K^{\rm T}B^{\rm T})$. It leads us to $B^{\rm T}P(PBK)^{\dagger}PBK = B^{\rm T}P$, which indicates that $R_{0} = -B^{\rm T}P(PBK)^{\dagger}PB$ is a particular solution. Then the general solution can be expressed by:
\begin{equation}\label{equ:positive_definite_R_Yr}
R_{1} = R_{0} + Y_{r} = -B^{\rm T}P(PBK)^{\dagger}PB + Y_{r},
\end{equation}
where $Y_{r} = Y_{r}^{\rm T}$ and $K^{\rm T}Y_{r} = 0$.

Secondly, we show that the matrix $R_{1}$ in the form of (\ref{equ:positive_definite_R_Yr}) is positive definite if and only if $Y_{r} = WW^{\rm T}$, where $W$ has $m - r_{b}$ linearly independent columns and $K^{\rm T}W = 0$. The analysis is broke into two cases. When $r_{b} = m$, the matrix $K^{\rm T}$ has full column rank. Then $Y_{r}$ must be a zero matrix and $R_{0}$ is the unique solution of $B^{\rm T}P = -RK$. For any non-zero vector $x_{2} \in \mathbb{R}^{m}$, it can be expressed as $x_{2} = K\rho$, where $\rho \in \mathbb{R}^{n}$ is non-zero. Then it holds that
\begin{equation*}
\begin{split}
x_{2}^{\rm T}R_{0}x_{2} & = -\rho^{\rm T}K^{\rm T}B^{\rm T}P(PBK)^{\dagger}PBK\rho = - \rho^{\rm T}PBK\rho.    
\end{split}
\end{equation*}
Since $-PBK \in \mathbb{S}^{n}_{+}$, it can be decomposed as $-PBK = Z^{\rm T}Z$ for some real $Z$, which leads to $BK\rho = -P^{-1}Z^{\rm T}Z\rho$. Since $K\rho \neq 0$ and $\text{rank}(BK) = \text{rank}(K)$, then $BK\rho \neq 0$ by rank-nullity theorem, which indicates that $Z^{\rm T}Z\rho \neq 0$. Then $x_{2}^{\rm T}R_{0}x_{2} = \rho^{\rm T}Z^{\rm T}Z\rho = \|Z\rho\|_{2}^{2} > 0$, which shows that $R_{0}$ is positive definite.

When $r_{b} < m$, to show its necessity, one can construct a vector as $x_{3} = (I - K(PBK)^{\dagger}PB)\eta$, where $\eta \in \mathbb{R}^{m}$. Since $R_{1}$ is positive definite, it holds that:
\begin{equation*}
\begin{split}
 x_{3}^{\rm T}R_{1}x_{3} & = x_{3}^{\rm T}R_{0}x_{3} + x_{3}^{\rm T}Y_{r}x_{3} \\ & = \eta^{\rm T}(I - B^{\rm T}P^{\rm T}(PBK)^{\dagger}K^{\rm T})Y_{r}(I - K(PBK)^{\dagger}PB) \eta \\ & = \eta^{\rm T}Y_{r}\eta \geq 0,
\end{split}
\end{equation*}
where the equality holds only for the case when $x_{3}$ is a zero vector. Then $Y_{r}$ is necessary to be positive semi-definite, which can be decomposed as $Y_{r} = WW^{\rm T}$. Furthermore, for any non-zero vector $\eta \in \nullspace(K^{\rm T})$, we denote $\eta^{\prime} = -K(PBK)^{\dagger}PB \eta$. Then $\eta^{\prime} \in \text{range}(K)$ and the vector $x_{3} = \eta + \eta^{\prime}$ is non-zero since $\eta$ and $\eta^{\prime}$ are orthogonal, which implies that $x_{3}^{\rm T}R_{1}x_{3} = \eta^{\rm T}Y_{r}\eta > 0$. We assume that the columns of $W$ do not span $\nullspace(K^{\rm T})$, then there exists a non-zero $\eta \in \nullspace(K^{\rm T})$ such that $x_{3}^{\rm T}R_{1}x_{3} = \|W^{\rm T}\eta\|_{2}^{2} = 0$, which is the contradiction.

To prove its sufficiency, it is equivalent to showing that for any non-zero vector $x_{4} \in \mathbb{R}^{m}$, it holds that $x_{4}^{\rm T}R_{1}x_{4} > 0$ for $R_{1}$ in (\ref{equ:matrixR_Ppositivedefinite}). By Lemma \ref{lemma:space_direct_sum}, the vector $x_{4}$ can be written as $x_{4} = \eta + K\rho$, where $\eta \in \nullspace(K^{\rm T})$. When $\eta = 0$, this can be shown by following the same steps as the proof for the case where $r_{b} = m$. When $\eta \neq 0$, it holds that
\begin{equation*}
x_{4}^{\rm T}Rx_{4} = -x_{4}^{\rm T}B^{\rm T}P(PBK)^{\dagger}PBx_{4} + \eta^{\rm T}WW^{\rm T}\eta.
\end{equation*}
Since $-PBK \in \mathbb{S}^{n}_{+}$, the first term is non-negative. It also holds that $\eta^{\rm T}WW^{\rm T}\eta = \|W^{\rm T}\eta\|_{2}^{2} > 0$ since the columns of $W$ span $\nullspace(K^{\rm T})$ and $\eta \in \nullspace(K^{\rm T})$. Then $x_{4}^{\rm T}R_{1}x_{4} > 0$, which proves the sufficiency.  \rulex

By Theorem \ref{theorem:givenP_findR}, a sufficient condition to uniquely recover $R$ by a single time point of $K(t)$ is given in Proposition \ref{proposition:GivenFQ_wellposedness}.


\begin{proposition}\label{proposition:GivenFQ_wellposedness}
Given $K(t)$ and $Q, F \in \mathbb{S}^{n}_{+}$, a sufficient condition for the well-posedness of Problem \ref{problem:1a} is that there exists $t_{1} \in [t_{0}, t_{f}]$ such that $K(t_{1})$ has full rank, and the unique solution of the matrix $R$ is given as:
\begin{equation}\label{equ:unique_Rmatrix_Kt_fullrank}
R_{0} = -B^{\rm T}P(t_{1})(P(t_{1})BK(t_{1}))^{\dagger}P(t_{1})B,
\end{equation}
which is the true value of $R$.
\end{proposition}

\proof Denote $R^{\star}$ as the true value of the control cost matrix. When there exists $t_{1} \in [t_{0}, t_{f}]$ such that $K(t_{1})$ has full column rank, the matrix $B$ has full column rank since $r_{b} \geq \rank(K(t_{1})B) = \rank(K(t_{1}))$. Then by Lemma \ref{lemma:necessary_sufficient_conditions_1973}, the matrix $P(t_{1})$ is positive definite. According to (\ref{equ:matrixR_Ppositivedefinite}), all the positive definite matrix $R(t_{1})$ satisfying $B^{\rm T}P(t_{1}) = -R(t_{1})K(t_{1})$ is given as
\begin{equation*}
R(t_{1}) = -B^{\rm T}P(t_{1})(P(t_{1})BK(t_{1}))^{\dagger}P(t_{1})B. 
\end{equation*}
Since the following equation is valid:
\begin{equation*}
\begin{split}
K(t_{1})^{\rm T}R(t_{1})K(t_{1}) = & -K(t_{1})^{\rm T}B^{\rm T}P(t_{1})(P(t_{1})BK(t_{1}))^{\dagger}P(t_{1})BK(t_{1}) \\ = & -P(t_{1})BK(t_{1}) = K(t_{1})^{\rm T}R^{\star}K(t_{1}),    
\end{split}
\end{equation*}
one can obtain that $K(t_{1})^{\rm T}(R(t_{1})-R^{\star})K(t_{1}) = 0$. Then $R(t_{1}) = R^{\star}$ follows from the assumption that $K(t_{1})$ has full rank. Thus, the true value $R^{\star}$ is recovered by using a single time point of $K(t)$.  \rulex

\begin{remark}\label{remark:R_time_variant_constant}
If there does not exist $t_{1}$ such that $K(t_{1})$ has full rank and $L_{1}$ is non-singular, one can always select a sequence of time points $\{t_{i}\}_{i=1}^{s} \subseteq [t_{0}, t_{f}]$ such that $\hat{K} = \begin{bmatrix} K_{t_{1}}, \ldots, K_{t_{s}}\end{bmatrix}$ has full row rank. Then the matrix $R$ can be unique recovered as $R = -(\hat{K}\hat{K}^{\rm T})^{-1}\hat{K}\hat{P}B$, where $\hat{P}^{\rm T} = \begin{bmatrix} P_{t_{1}}^{\rm T}, \ldots, P_{t_{s}}^{\rm T} \end{bmatrix}$. Comparing to the first method, the second method enables us to evaluate $R$ using some time points of $K(t)$, or even a single time point, rather than relying on the full trajectory of $K(t)$ over $[t_{0}, t_{f}]$. This significantly simplifies data requirements and reduces computational burden. Particularly, in scenarios where the data is collected in real time, estimation of the matrix $R$ can also be updated in real time with reduced uncertainty, as soon as  $\hat{K}$ is updated.
\end{remark}

What remains to be discussed is the case where $P(\tilde{t})$ is positive semi-definite for $\tilde{t} \in [t_{0}, t_{f}]$. According to Lemma 
\ref{lemma:necessary_sufficient_conditions_1973}, the matrix $K(\tilde{t})$ is rank deficient. In Proposition \ref{proposition:solutionR_Ppositivesemi-definite}, the solution space of $R(\tilde{t})$ will be investigated.

\begin{proposition}\label{proposition:solutionR_Ppositivesemi-definite}
For $\tilde{t} \in [t_{0}, t_{f}]$ such that $P(\tilde{t})$ is a positive semi-definite matrix, all the symmetric solution of $B^{\rm T}P(\tilde{t}) = -R(\tilde{t})K(\tilde{t})$ are represented in terms of $P(\tilde{t})$ by $R(\tilde{t}) = R_{2}(\tilde{t}) +Y$. Omitting $'(\tilde{t})'$ for brevity, the matrix $R_{2}$ takes the form
\begin{equation}\label{equ:particular_R1_K_dagger}
R_{2} = -B^{\rm T}PK^{\dagger} - (K^{\dagger})^{\rm T}PB + (K^{\dagger})^{\rm T}K^{\rm T}B^{\rm T}PK^{\dagger}, 
\end{equation}
and $Y = Y^{\rm T}$, $K^{\rm T}Y = 0$. Then the solution set of $R(\tilde{t})$, denoted by $\mathcal{V}_{R}(\tilde{t})$, is given as below:
\begin{equation}\label{equ:solution_set_R_positivesemi-definite_P}
\mathcal{V}_{R}(\tilde{t}) = \{R \mid R = R_{2}(\tilde{t}) + Y(\tilde{t})\} \cap \mathbb{S}^{m}_{++},
\end{equation}
to further ensure the positive definiteness of $R$.
\end{proposition}

\proof See Appendix \ref{appendix:R_P_posisemi}.  \rulex

After deriving $R$ based on $Q$ and $F$, we will further study Problem \ref{problem:1b}, where only $F$ is given, while $Q$ remains unknown. Since $P(t)$ can not be computed explicitly, the first method is not applicable, whereas the second method can be naturally applied. In the rest of this paper, we denote $K(t_{f})$ as $K_{f}$. Consider $B^{\rm T}P(t) = -RK(t)$ on $t = t_{f}$, it holds that $B^{\rm T}F = -RK_{f}$. Then the solution space of $R$ is discussed in Proposition \ref{proposition:GivenFstar_uniqueRstar}.


\begin{proposition}\label{proposition:GivenFstar_uniqueRstar}
Given $K(t)$ and $F \in \mathbb{S}^{n}_{+}$, if it satisfies $\rank(FB) = m$, then $R$ can be uniquely recovered as:
\begin{equation}\label{equ:Rsolutionspace_F_posi}
R_{0} = -B^{\rm T}F(FBK_{f})^{\dagger}FB.
\end{equation}
Furthermore, the equality $\rank(FB) = m$ holds if and only if $B$ has full column rank and $F \in \mathbb{S}^{n}_{++}$.

If $\rank(FB) < m$ and $\rank(K_{f}) = r_{b}$, the solution space of $R$ is explicitly characterized as:
\begin{equation}\label{equ:solution_space_R_F_psd}
R_{1} = -B^{\rm T}F(FBK_{f})^{\dagger}FB + WW^{\rm T},    
\end{equation}
where $W$ has $m-r_{b}$ linearly independent columns and $K_{f}^{\rm T}W = 0$. If $\rank(FB) < m$ and $\rank(K_{f}) < r_{b}$, the solution space of $R$ can be evaluated as $\mathcal{V}_{R}(t_{f})$ in (\ref{equ:solution_set_R_positivesemi-definite_P}).
\end{proposition}

\proof Since $\rank(FB) = m$, it can be obtained that $\rank(K_{f}^{\rm T}R) = m$, which implies that $\rank(K_{f}) = m$ since $R$ must be positive definite. With $K_{f}$ having full rank, by Proposition  \ref{proposition:GivenFQ_wellposedness}, the positive definite matrix $R$ to satisfy $B^{\rm T}F = -RK_{f}$ can be uniquely given by (\ref{equ:Rsolutionspace_F_posi}).

Next, we derive the necessary and sufficient condition for $\rank(FB) = m$. When $B$ has full column rank and $F \in \mathbb{S}^{n}_{++}$, it holds that $\rank(FB) = \rank(B) = m$, which shows the sufficiency. Then we prove its necessity. If $B$ does not have full rank, it holds that $\rank(FB) \leq r_{b} < m$. On the other hand, if $F$ is only positive semi-definite, then $\rank(FB) = \rank(K_{f}) < r_{b} \leq m$ by Lemma \ref{lemma:necessary_sufficient_conditions_1973}.

When $\rank(FB) < m$, it is clear that $K_{f}$ is rank deficient. If $\rank(K_{f}) = r_{b}$, the matrix $F$ is positive definite by Lemma \ref{lemma:necessary_sufficient_conditions_1973}, then the solution space of $R$ shown in (\ref{equ:solution_space_R_F_psd}) can be derived by applying Theorem \ref{theorem:givenP_findR} at $\tilde{t} = t_{f}$. If $\rank(K_{f}) < r_{b}$, the matrix $F$ is only positive semi-definite, then the solution space $\mathcal{V}_{R}(t_{f})$ can be characterized by applying Proposition \ref{proposition:solutionR_Ppositivesemi-definite} at $\tilde{t} = t_{f}$.  \rulex
%


\section{Reconstruction of $Q$ and $F$}\label{sec:Rknown}

In this section, we investigate the reconstruction of $Q$ and $F$ when $R$ is given in a more general setup. Following \citep{li2020continuous}, which addresses the existence problem under the assumption of system controllability, we introduce several key extensions in this work. Firstly, we clarify the role of system controllability in establishing the bijective mapping between $Q$ and $F$, i.e., the well-posedness of IOC problems. Secondly, Problem \ref{problem:2a} is studied. When the bijectivity condition holds, the analytical expression is provided; otherwise, the solution space is characterized via LMI formulations. Thirdly, for the cases where $Q$ and $F$ are both unknown, i.e., Problem \ref{problem:2b}, an explicit sufficient condition for well-posedness is provided, along with analytical expressions for $Q$ and $F$.

To study the mapping between $Q$ and $F$, we assume that $(Q_{1}, F_{1})$ and $(Q_{2}, F_{2})$ are two pairs of solutions. The corresponding solution to (DRE) are denoted by $P_{1}(t)$ and $P_{2}(t)$, respectively, i.e., for $t \in [t_{0}, t_{f}]$, it holds that
\begin{equation*}
\begin{cases}
-\dot{P}_{1} = A^{\rm T}P_{1} + P_{1}A + P_{1}BK + Q_{1}, \, P_{1}(t_{f}) = F_{1}, \\
-\dot{P}_{2} = A^{\rm T}P_{2} + P_{2}A + P_{2}BK + Q_{2}, \, P_{2}(t_{f}) = F_{2}, \\
-RK(t) = B^{\rm T}P_{1}(t) = B^{\rm T}P_{2}(t).
\end{cases}
\end{equation*}
Denote ${\it \Delta} P(t) = P_{1}(t) - P_{2}(t)$, ${\it \Delta} Q = Q_{1} - Q_{2}$ and ${\it \Delta} F = F_{1} - F_{2}$. By vectorization, it holds that for $t \in [t_{0},t_{f}]$:
\begin{equation}\label{equ:vectorization_givenR_deltaP_deltaQ}
\begin{split}
\vectorize({\it \Delta} \dot{P}(t)) & = \tilde{A}^{\rm T} \vectorize({\it \Delta} P(t)) - \vectorize({\it \Delta} Q), \\ \tilde{B}^{\rm T}\vectorize({\it \Delta} P(t)) & = 0, \\ \vectorize({\it \Delta} P(t_{f})) & = \vectorize({\it \Delta} F),
\end{split}
\end{equation}
where the matrices $\tilde{A}$ and $\tilde{B}$ are given as:
\begin{equation}\label{equ:matrix_tildeA_tildeB}
\begin{split}
\tilde{A} = -I_{n} \otimes A - A \otimes I_{n}, \quad  \tilde{B} = I_{n} \otimes B.      
\end{split}
\end{equation}
By taking $n^{2}$ derivatives of $\tilde{B}^{\rm T}\vectorize({\it \Delta} P(t)) = 0$, one can obtain that for all $t \in [t_{0},t_{f}]$, it holds that
\begin{equation}\label{equ:vectorization_deltaP_deltaQ}
\tilde{H} \vectorize({\it \Delta} P(t)) + \tilde{N} \vectorize({\it \Delta} Q) = 0,
\end{equation}
where $\vectorize({\it \Delta} P(t_{f})) = \vectorize({\it \Delta} F)$ and $\tilde{H}$, $\tilde{N}$ are defined as: 
\begin{equation}\label{equ:tilde_H_tilde_N}
\begin{split}
\tilde{H} & = \begin{bmatrix} \tilde{B}, \tilde{A}\tilde{B}, \ldots, \tilde{A}^{n^{2}-1}\tilde{B}, \tilde{A}^{n^{2}}\tilde{B}\end{bmatrix}^{\rm T}, \\
\tilde{N} & = \begin{bmatrix}
0, -\tilde{B}, -\tilde{A}\tilde{B}, \ldots, -\tilde{A}^{n^{2}-1}\tilde{B} \end{bmatrix}^{\rm T}.
\end{split}
\end{equation}
To gain further insight into $\tilde{H}$ and $\tilde{N}$, a new linear system is constructed as $\dot{\tilde{x}} =  \tilde{A}\tilde{x} + \tilde{B}\tilde{u}$,
where $\tilde{A}$ and $\tilde{B}$ are given in (\ref{equ:matrix_tildeA_tildeB}). Its controllability matrix is defined as:
\begin{equation}\label{equ:controllability_matrix_constructed_system}
\tilde{{\it \Gamma}}_{c} = \begin{bmatrix}\tilde{B}, \tilde{A}\tilde{B}, \ldots, \tilde{A}^{n^{2}-1}\tilde{B}\end{bmatrix}.
\end{equation}

The Cayley-Hamilton theorem \citep{petersen2012linearslgebra} will be used in the following proofs, which is introduced first in the next lemma.

\begin{lemma}\label{lemma:cayley_hamilton}
Let $A \in \mathbb{R}^{n \times n}$. The matrix $A^{k}$, where $k \geq n$, can be expressed as the linear combination of the matrices $I$, $A$, $A^{2}$, $\ldots$ and $A^{n-1}$.
\end{lemma}

The matrices $\tilde{H}$ and $\tilde{N}$ can be written as $\tilde{H} = \begin{bmatrix} \tilde{{\it \Gamma}}_{c}, \, \tilde{A}^{n^{2}}\tilde{B}\end{bmatrix}^{\rm T}$ and $\tilde{N} = \begin{bmatrix} 0,\, -\tilde{{\it {\it \Gamma}}}_{c}\end{bmatrix}^{\rm T}$. It is clear that $\nullspace(\tilde{N}) = \nullspace(\tilde{{\it {\it \Gamma}}}_{c}^{\rm T})$, and $\nullspace(\tilde{H}) = \nullspace(\tilde{{\it {\it \Gamma}}}_{c}^{\rm T})$ by Cayley-Hamilton theorem. Then the matrices $\tilde{H}$, $\tilde{N}$ and $\tilde{{\it {\it \Gamma}}}_{c}^{\rm T}$ share the same null space. Furthermore, the properties of $\tilde{{\it {\it \Gamma}}}_{c}$ are established in the next proposition. In the rest of the paper, we denote $D_{n}$ as the duplication matrix to transform half-vectorization into vectorization for $n\times n$ symmetric matrices.

\begin{proposition}\label{Proposition:tilde_Gamma_controllability}
The following statements about $\tilde{{\it {\it \Gamma}}}_{c}$ in (\ref{equ:controllability_matrix_constructed_system}) are valid:
\begin{enumerate}
\item  The matrix $\tilde{{\it {\it \Gamma}}}_{c}^{\rm T}$ has full column rank if and only if $(A, B)$ is controllable,
\item The matrix $\tilde{{\it {\it \Gamma}}}_{c}^{\rm T}D_{n}$ has full column rank if and only if $(A, B)$ is controllable.
\end{enumerate}
\end{proposition}

\proof We begin with the first statement. According to \citep{li2020continuous}, if $(A, B)$ is controllable, the matrix $\begin{bmatrix} \tilde{B}, \tilde{A}\tilde{B}, \ldots, \tilde{A}^{n-1}\tilde{B} \end{bmatrix}^{\rm T}$ has full column rank. It indicates that $\tilde{{\it {\it \Gamma}}}_{c}^{\rm T}$ also has full column rank, which shows the sufficiency. Its necessity is established by proving the contrapositive statement.  If $(A,B)$ is not controllable, by PBH test, there exists a non-zero vector $w^{\rm T}$ such that $w^{\rm T}A = \lambda w^{\rm T}$ and $w^{\rm T}B = 0$. Denote $\tilde{w}^{\rm T} = w^{\rm T} \otimes w^{\rm T}$. The following equations are satisfied:
\begin{equation*}
\begin{split}
\tilde{w}^{\rm T}\tilde{A} & = -(w^{\rm T} \otimes w^{\rm T})(I_{n} \otimes A  + A \otimes I_{n}) = - w^{\rm T} \otimes w^{\rm T}A - w^{\rm T}A \otimes w^{\rm T} = -2\lambda \tilde{w}^{\rm T}, \\
\tilde{w}^{\rm T}\tilde{B} & = (w^{\rm T} \otimes w^{\rm T})(I_{n} \otimes B) = w^{\rm T} \otimes 0 = 0,
\end{split}
\end{equation*}
which indicates that $(\tilde{A},\tilde{B})$ is not controllable. Then the matrix $\tilde{{\it {\it \Gamma}}}_{c}^{\rm T}$ in (\ref{equ:controllability_matrix_constructed_system}) does not have full column rank.

Then we proceed to the second statement. If $(A, B)$ is controllable, the matrix $\tilde{{\it {\it \Gamma}}}_{c}^{\rm T}$ has full column rank, which indicates that $\rank(\tilde{{\it {\it \Gamma}}}^{\rm T}_{c}D_{n}) = \rank(D_{n})$. Since the duplication matrix $D_{n}$ has full column rank, it implies that $\tilde{{\it \Gamma}}^{\rm T}_{c}D_{n}$ also has full column rank. Furthermore, if $(A, B)$ is not controllable, there exists a non-zero vector $v \in \mathbb{R}^{n}$ such that $\begin{bmatrix} B, AB, \ldots, A^{n-1}B \end{bmatrix}^{\rm T}v = 0$. By Cayley-Hamilton theorem, the vector $v$ satisfies $B^{\rm T}(A^{i})^{\rm T}v = 0$ for all positive integer $i$. For any integer $0 \leq k \leq n^{2}-1$, the term $\tilde{A}^{k}\tilde{B}$ can be expressed as 
\begin{equation*}
\begin{split}
\tilde{A}^{k}\tilde{B} & = (-1)^{k}(I_{n} \otimes A + A \otimes I_{n})^{k}(I_{n} \otimes B) \\ & = (-1)^{k}(\binom{k}{0}I_{n} \otimes A^{k} + \binom{k}{1}A \otimes A^{k-1} + \ldots + \binom{k}{k}A^{k} \otimes I_{n})(I_{n} \otimes B) \\ & = (-1)^{k}(\binom{k}{0}I_{n} \otimes A^{k}B + \binom{k}{1}A \otimes A^{k-1}B + \ldots + \binom{k}{k}A^{k} \otimes B) \\ & = (-1)^{k} \sum_{j = 0}^{k} \binom{k}{j} A^{j} \otimes A^{k-j}B.
\end{split}    
\end{equation*}
Then the following equation holds for all integer $0 \leq k \leq n^{2}-1$:
\begin{equation}
\begin{split}
\tilde{B}^{\rm T} (\tilde{A}^{\rm T})^{k} \vectorize(vv^{\rm T}) = & (-1)^{k} \sum_{j = 0}^{k} \binom{k}{j} (A^{\rm T})^{j} \otimes B^{\rm T}(A^{\rm T})^{k-j}\vectorize(vv^{\rm T}) \\ = & (-1)^{k} \sum_{j = 0}^{k} \binom{k}{j} \vectorize(B^{\rm T}(A^{\rm T})^{k-j}vv^{\rm T}A^{j}) = 0.
\end{split}
\end{equation}
Since $v$ is a non-zero vector, then $vv^{\rm T}$ is a matrix with rank $1$, then $\halfvectorize(vv^{\rm T})$ is also non-zero. Then it holds that:
\begin{equation}
\begin{split}
\tilde{{\it {\it \Gamma}}}_{c}^{\rm T}D_{n} \halfvectorize(vv^{\rm T}) = \tilde{{\it {\it \Gamma}}}_{c}^{\rm T} \vectorize(vv^{\rm T}) = \begin{bmatrix} \tilde{B}, \tilde{A}\tilde{B}, \ldots, \tilde{A}^{n^{2}-1}\tilde{B} 
\end{bmatrix}^{\rm T}\vectorize(vv^{\rm T}) = 0,
\end{split}
\end{equation}
which proves that $\tilde{{\it {\it \Gamma}}}_{c}^{\rm T}D_{n}$ does not have full column rank.  \rulex

To investigate the mapping between $Q$ and $F$ with the given $R$, we first analyse the solution space of $Q$ (or $F$) with a fixed $F$ (or $Q$) in the next theorem.

\begin{theorem}\label{theorem:GivenR_QF_solution_space}
Given $K(t)$ and $R \in \mathcal{R}_{s}$, let $Q_{1}$ and $F_{1}$ denote a pair of solution. With $F_{1}$ fixed, the solution space of $Q$ is $\mathcal{D}_{Q+} = \mathcal{D}_{Q} \cap S_{+}^{n}$, where $\mathcal{D}_{Q}$ is defined as  
\begin{equation}\label{equ:GivenR_solution_space_Q}
\begin{split}
\mathcal{D}_{Q} = \{Q \mid Q = Q_{1} + {\it \Delta} Q, \, \tilde{{\it {\it \Gamma}}}_{c}^{\rm T}\vectorize({\it \Delta} Q) = 0\}. 
\end{split}
\end{equation}
Similarly, with $Q_{1}$ fixed, the solution space of $F$ is $\mathcal{D}_{F+} = \mathcal{D}_{F} \cap S_{+}^{n}$, where $\mathcal{D}_{F}$ is defined as
\begin{equation}\label{equ:GivenR_solution_space_F}
\begin{split}
\mathcal{D}_{F} = \{F \mid F = F_{1} + {\it \Delta} F, \, \, \tilde{{\it {\it \Gamma}}}_{c}^{\rm T}\vectorize({\it \Delta} F) = 0\},
\end{split}
\end{equation}
where $\tilde{{\it {\it \Gamma}}}_{c}$ is shown in (\ref{equ:controllability_matrix_constructed_system}).
\end{theorem}


\proof Denote $P_{1}(t)$ as the solution of (DRE) corresponding with $Q_{1}$ and $F_{1}$. Let $Q_{2}$ and $F_{2}$ denote another pair of solution and $P_{2}(t)$ denote the solution of the corresponding (DRE). By (\ref{equ:vectorization_givenR_deltaP_deltaQ}), one can express $\vectorize({\it \Delta} P(t))$, where $t \in [t_{0}, t_{f}]$, by $\vectorize({\it \Delta} Q)$ and $\vectorize({\it \Delta} F)$ as below:
\begin{equation}\label{equ:vecDeltaP_expression}
\vectorize({\it \Delta} P(t)) =  {\rm e}^{\tilde{A}^{\rm T}(t-t_{f})}\vectorize({\it \Delta} F)  - \int_{t_{f}}^{t} \! {\rm e}^{\tilde{A}^{\rm T}(t-s)}\vectorize({\it \Delta} Q)  \, \mathrm{d}s.
\end{equation}

Firstly, we show that the solution space of $Q$ is $\mathcal{D}_{Q+}$ when ${\it \Delta} F = 0$. The necessity is proved by considering (\ref{equ:vectorization_deltaP_deltaQ}) at $t = t_{f}$. Then $\vectorize({\it \Delta} Q)$ must satisfy $\tilde{N}\vectorize({\it \Delta} Q) = 0$, i.e.,  $\vectorize({\it \Delta} Q) \in \nullspace(\tilde{{\it {\it \Gamma}}}_{c}^{\rm T})$. Next, we prove its sufficiency. By substituting (\ref{equ:vecDeltaP_expression}) into $\tilde{B}^{\rm T}\vectorize(P(t)) = 0$ and plugging in ${\it \Delta} F = 0$, one can obtain that if $\vectorize({\it \Delta} Q) \in \nullspace(\tilde{{\it {\it \Gamma}}}_{c}^{\rm T})$, it holds that:
\begin{equation*}
\begin{split}
\tilde{B}^{\rm T} \vectorize({\it \Delta} P(t)) =  -\tilde{B}^{\rm T}\int_{t_{f}}^{t} \! {\rm e}^{\tilde{A}^{\rm T}(t-s)}\vectorize({\it \Delta} Q)  \, \mathrm{d}s  =    -\sum_{k = 0}^{\infty} \int_{t_{f}}^{t} \! \tfrac{(t-s)^{k}}{k!} \tilde{B}^{\rm T}(\tilde{A}^{\rm T})^{k}\vectorize({\it \Delta} Q)   \, \mathrm{d}s = 0,
\end{split}
\end{equation*}
since $\tilde{B}^{\rm T}(\tilde{A}^{\rm T})^{k} \vectorize({\it \Delta} Q) = 0$ for all $k \geq n^{2}$ by Cayley–Hamilton theorem. When $F_{1}$ is fixed, the corresponding solution of $Q$ must be in the form of $Q = Q_{1} + {\it \Delta} Q$, where $\vectorize({\it \Delta} Q) \in \nullspace(\tilde{{\it {\it \Gamma}}}_{c}^{\rm T})$. To further ensure positive semi-definiteness of $Q_{1} + {\it \Delta} Q$, the solution space of $Q$ is $\mathcal{D}_{Q+} = \mathcal{D}_{Q} \cap S_{+}^{n}$, where $\mathcal{D}_{Q}$ is shown in (\ref{equ:GivenR_solution_space_Q}).


Secondly, when ${\it \Delta} Q = 0$, the solution space of $F$ is shown to be $\mathcal{D}_{F+}$. Consider (\ref{equ:vectorization_deltaP_deltaQ}) at $t = t_{f}$, one can obtain that $\vectorize({\it \Delta} F)$ must satisfy $\tilde{H}\vectorize({\it \Delta} F) = 0$, i.e., $\vectorize({\it \Delta} F) \in \nullspace(\tilde{{\it {\it \Gamma}}}_{c}^{\rm T})$, which proves the necessity. Its sufficiency is proved by substituting (\ref{equ:vecDeltaP_expression}) into $\tilde{B}^{\rm T}\vectorize(P(t)) = 0$ and plugging in ${\it \Delta} Q = 0$, by Cayley–Hamilton theorem, if it satisfies $\vectorize({\it \Delta} F) \in \nullspace(\tilde{{\it {\it \Gamma}}}^{\rm T}_{c})$, it holds that:
\begin{equation*}
\begin{split}
\tilde{B}\vectorize({\it \Delta} P(t)) =  \tilde{B}^{\rm T}{\rm e}^{\tilde{A}^{\rm T}(t-t_{f})}\vectorize({\it \Delta} F)  = \sum_{k = 0}^{\infty} \tfrac{(t-t_{f})^{k}}{k!} \tilde{B}^{\rm T}(\tilde{A}^{\rm T})^{k}  \vectorize({\it \Delta} F) = 0.
\end{split}
\end{equation*}
Therefore, when the matrix $Q_{1}$ is fixed, the corresponding solution of $F$ must be in the form of $F = F_{1} + {\it \Delta} F$, where $\vectorize({\it \Delta} F) \in \nullspace(\tilde{{\it {\it \Gamma}}}_{c}^{\rm T})$. Then the solution space of $F$ is given as $\mathcal{D}_{F+} = \mathcal{D}_{F} \cap S_{+}^{n}$, where $\mathcal{D}_{F}$ is shown in (\ref{equ:GivenR_solution_space_F}), to guarantee positive semi-definiteness.  \rulex

Since the matrices $Q$ and $F$ should be symmetric, the vectorization of $Q$ and $F$ can be expressed as $\vectorize({\it \Delta} Q) = D_{n}\halfvectorize({\it \Delta} Q)$ and $\vectorize({\it \Delta} F) = D_{n}\halfvectorize({\it \Delta} F)$. Previous studies have established that the set $\mathcal{D}_{Q}$ (or $\mathcal{D}_{F}$) only contains a single element if the underlying system $(A,B)$ is controllable. By Proposition \ref{Proposition:tilde_Gamma_controllability}, we further show that the controllability of $(A, B)$ is also the necessary condition. However, whether it can serve as the necessary condition for the bijective mapping between $Q$ and $F$, that is, whether $\mathcal{D}_{Q+}$ (or $\mathcal{D}_{F+}$) only contains a single element, depends on the specific properties of the given $Q_{1}$ (or $F_{1}$) and the null space of $\tilde{{\it {\it \Gamma}}}_{c}^{\rm T}$. For example, if $Q_{1}$ is positive definite, the controllability of $(A, B)$ is a necessary condition; in contrast, if $Q_{1}$ is given as the zero matrix and ${\it \Delta} Q$ that satisfies $\tilde{{\it {\it \Gamma}}}_{c}^{\rm T}D_{n}\halfvectorize({\it \Delta} Q) = 0$ is indefinite, the controllability is not necessary.

After investigating the mapping between $Q$ and $F$, we proceed to explore the solution space for $Q$ (or $F$) when $R$ and $F$ (or $Q$) are given. In \citep{li2020continuous}, (DRE) corresponding with (\ref{equ:classic_optimal_problem}) is reformulated into a Lyapunov equation, which is summarized in Lemma \ref{lemma:Pt_representation_by_R}.

\begin{lemma}\label{lemma:Pt_representation_by_R}
Denote $P_{0}(t)$ and $G(t)$ as:
\begin{equation*}
\begin{split}
& P_{0}(t)  = -K(t)^{\rm T}R(RK(t)B)^{\dagger}RK(t), \\
& G(t)  = A^{\rm T}P_{0}(t) + P_{0}(t)A + P_{0}(t)BK(t) + \dot{P}_{0}(t).
\end{split}
\end{equation*}
The solution $P(t)$ of (DRE) with the given $R$ can be expressed as $P(t) = P_{0}(t) + Y(t)$, where $Y(t)$ satisfies: 
\begin{equation}\label{equ:GivenR_vectorYt}
\dot{Y}(t)  = -A^{\rm T}Y(t) - Y(t)A - Q - G(t),
\end{equation}
with $Y(t_{f}) = F - P_{0}(t_{f})$ and $B^{\rm T}Y(t) = 0$.
\end{lemma}

We denote $F_{0} = P_{0}(t_{f}) =  -K_{f}^{\rm T}R(RK_{f}B)^{\dagger}RK_{f}$. In the next theorem, we first discuss the case where $(A, B)$ is controllable.


\begin{theorem}\label{theorem:GivenRstar_QFstar_unique}
Let $g(t)^{\rm T} = \begin{bmatrix} g_{0}^{\rm T}, g_{1}^{\rm T}(t), \ldots, g_{n^{2}}^{\rm T}(t)\end{bmatrix}$, where $g_{0} \in \mathbb{R}^{nm \times 1}$ is a zero vector and $g_{k}(t) \in \mathbb{R}^{nm \times 1}$, for $k \in [n^{2}]$, is defined as
\begin{equation*}
g_{k}(t) = -\sum_{i=0}^{k-1}\tilde{B}^{\rm T}(\tilde{A}^{\rm T})^{i} \vectorize(G(t)^{(k-1-i)}).
\end{equation*} 
When $(A, B)$ is controllable, given $K(t)$ and $R \in \mathcal{R}_{s}$, if $F$ is further given, the matrix $Q$ can be uniquely recovered as:
\begin{equation}\label{equ:Qstar_RFgiven_ABcontrol}
Q  = -\matrixize(D_{n}(\tilde{N}^{\rm T}_{d_{n}}\tilde{N}_{d_{n}})^{-1}\tilde{N}^{\rm T}_{d_{n}}[\tilde{H}_{d_{n}}y_{f} + g(t_{f})]),
\end{equation}
where $y_{f} = \halfvectorize(F - F_{0})$.

Similarly, if $Q$ is further given, the matrix $F$ can be uniquely recovered as $F = F_{0} + \tilde{F}$, where
\begin{equation}\label{equ:Fstar_RQgiven_ABcontrol}
\tilde{F} = -\matrixize(D_{n}(\tilde{H}_{d_{n}}^{\rm T}\tilde{H}_{d_{n}})^{-1}\tilde{H}_{d_{n}}^{\rm T}[\tilde{N}_{d_{n}}q + g(t_{f})]),
\end{equation}
with $q = \halfvectorize(Q)$.
\end{theorem}

\proof By vectorizing (\ref{equ:GivenR_vectorYt}), it holds for all $t \in [t_{0}, t_{f}]$:
\begin{equation*}
\begin{cases}
\vectorize{(\dot{Y}(t))}  = \tilde{A}^{\rm T}\vectorize(Y(t)) - \vectorize(Q) - \vectorize(G(t)), \\
\tilde{B}^{\rm T}\vectorize(Y(t)) = 0,
\end{cases}
\end{equation*}
with the terminal condition $\vectorize(Y(t_{f})) = \vectorize(F - F_{0})$. For any positive integer $k$, the $k$-th derivative of the equation $\tilde{B}^{\rm T}\vectorize(Y(t)) = 0$ is given as below:
\begin{equation*}
\begin{split}
0 =  \tilde{B}^{\rm T}(\tilde{A}^{\rm T})^{k}\vectorize(Y(t)) - \tilde{B}^{\rm T}(\tilde{A}^{\rm T})^{k-1}\vectorize(Q) - \sum_{i=0}^{k-1}\tilde{B}^{\rm T}(\tilde{A}^{\rm T})^{i}\vectorize(G(t)^{(k-1-i)}).   
\end{split}
\end{equation*}
Then by taking $n^{2}$ derivatives of $\tilde{B}^{\rm T}\vectorize(Y(t)) = 0$, the following equation holds for all $t \in [t_{0}, t_{f}]$:
\begin{equation}\label{equ:nsquare_derivatives_Yt}
\tilde{H}\vectorize(Y(t)) + \tilde{N}\vectorize(Q) + g(t) = 0.
\end{equation}
where $\tilde{H}$ and $\tilde{N}$ are shown in (\ref{equ:tilde_H_tilde_N}). Since both of $Y(t)$ and $Q$ are symmetric, then (\ref{equ:nsquare_derivatives_Yt}) can be reformulated as:
\begin{equation}\label{equ:nsquare_derivatives_Yt_dn}
\tilde{H}_{d_{n}} \halfvectorize(Y(t)) + \tilde{N}_{d_{n}}  \halfvectorize(Q) + g(t) = 0,
\end{equation}
where $\tilde{H}_{d_{n}} = \tilde{H}D_{n}$ and $\tilde{N}_{d_{n}} = \tilde{N}D_{n}$. Since $\tilde{N}_{d_{n}}\halfvectorize(Q)$ is constant, then $\tilde{H}_{d_{n}}\halfvectorize(Y(t)) + g(t)$ must also remain constant over $t \in [t_{0}, t_{f}]$. By considering (\ref{equ:nsquare_derivatives_Yt_dn}) on $t = t_{f}$, it holds that
\begin{equation}\label{equ:nsquare_derivatives_Yt_f}
\tilde{H}_{d_{n}} \halfvectorize(Y(t_{f})) + \tilde{N}_{d_{n}}  \halfvectorize(Q) + g(t_{f}) = 0.
\end{equation}
Since $D_{n}$ has full column rank, it is clear that $\tilde{H}_{d_{n}}$, $\tilde{N}_{d_{n}}$ and $\tilde{{\it \Gamma}}_{c}^{\rm T}D_{n}$ share the same null space, which means that $\tilde{H}_{d_{n}}$ and $\tilde{N}_{d_{n}}$ have full column rank.

If $F$ is further given, according to (\ref{equ:nsquare_derivatives_Yt_f}), it holds that:
\begin{equation}\label{equ:GivenR_Fstar_further_Q}
\tilde{N}_{d_{n}}\halfvectorize(Q) = - \tilde{H}_{d_{n}}y_{f} - g(t_{f}),
\end{equation}
where $y_{f} = \halfvectorize(F - F_{0})$. Then the unique solution is given as 
\begin{equation}
\halfvectorize(Q) = - (\tilde{N}^{\rm T}_{d_{n}}\tilde{N}_{d_{n}})^{-1}\tilde{N}^{\rm T}_{d_{n}}[\tilde{H}_{d_{n}}y_{f} + g(t_{f})].    
\end{equation}
The analytical expression of $Q$ in (\ref{equ:Qstar_RFgiven_ABcontrol}) follows.

Similarly, if $Q$ is further given, one can rewrite (\ref{equ:nsquare_derivatives_Yt_f}) as:
\begin{equation}\label{equ:GivenR_Qstar_further_F}
\tilde{H}_{d_{n}}\halfvectorize(F - F_{0}) = -\tilde{N}_{d_{n}}q - g(t_{f}),
\end{equation}
where $q = \halfvectorize(Q)$. Then it holds that $\halfvectorize(F-F_{0}) = -(\tilde{H}_{d_{n}}^{\rm T}\tilde{H}_{d_{n}})^{-1}\tilde{H}_{d_{n}}^{\rm T}[\tilde{N}_{d_{n}}q + g(t_{f})]$. The analytical expression of $F$ follows, which completes the proof. \rulex

When $(A, B)$ is not controllable, the matrix $Q$ (or $F$) may not be recovered uniquely when $R$ and $F$ (or $Q$) are given. In this case, its solution space will be characterized by using the tool of LMI in Proposition \ref{proposition:GivenRstar_QFstar_notunique}.


\begin{proposition}\label{proposition:GivenRstar_QFstar_notunique}
When $(A, B)$ is not controllable, let ${\it \Gamma}_{1}, \ldots, {\it \Gamma}_{r}$ span the null space of $\tilde{{\it \Gamma}}_{c}^{\rm T}D_{n}$. Given $K(t)$ and $R \in \mathcal{R}_{s}$, if $F$ is further given, the solution space of the matrix $Q$ is given as
\begin{equation}\label{equ:GivenRF_DQplus_not_unique}
Q = \matrixize(D_{n}q + \sum_{i=1}^{r}\alpha_{i}D_{n}{\it \Gamma}_{i}),
\end{equation}
where the vector $q$ satisfies 
\begin{equation}\label{equ:DQplus_not_unique_q_particular}
\tilde{N}_{d_{n}}q = - \tilde{H}_{d_{n}}\halfvectorize(F - F_{0}) - g(t_{f}),
\end{equation}
and $\alpha_{1}, \ldots, \alpha_{r}$  solve the following LMI problem:
\begin{equation}\label{equ:DQplus_not_unique_q_LMI}
\matrixize(D_{n}q) + \sum_{i = 1}^{r}\alpha_{i} \matrixize(D_{n}{\it \Gamma}_{i}) \in \mathbb{S}^{n}_{+}.
\end{equation}
If $Q$ is further given, the solution space of $F$ is given as
\begin{equation}\label{equ:GivenRQ_DFplus_not_unique}
F = F_{0} + \matrixize(D_{n}f + \sum_{i=1}^{r}\beta_{i}D_{n}{\it \Gamma}_{i}),
\end{equation}
where the vector $f$ satisfies
\begin{equation}\label{equ:DQplus_not_unique_f_particular}
\tilde{H}_{d_{n}}f = -\tilde{N}_{d_{n}} \halfvectorize(Q) - g(t_{f}),
\end{equation}
and $\beta_{1}, \ldots, \beta_{r}$ solve the following LMI problem:
\begin{equation}\label{equ:DQplus_not_unique_f_LMI}
F_{0} + \matrixize(D_{n}f) + \sum_{i = 1}^{r}\beta_{i}\matrixize(D_{n}{\it \Gamma}_{i})  \in \mathbb{S}^{n}_{+}.
\end{equation}
\vspace{-0.6cm}
\end{proposition}

\proof In this case, by following the same steps in the proof of Theorem \ref{theorem:GivenRstar_QFstar_unique}, one can obtain that (\ref{equ:nsquare_derivatives_Yt_f}) is still valid but $\tilde{H}_{d_{n}}$ and $\tilde{N}_{d_{n}}$ in (\ref{equ:nsquare_derivatives_Yt_f}) do not have full column rank. If $F$ is further specified, the matrix $Q$ satisfies (\ref{equ:GivenR_Fstar_further_Q}). Then $\halfvectorize(Q)$ must be of the form $\halfvectorize(Q) = q + \sum_{i=1}^{r} \alpha_{i}{\it \Gamma}_{i}$, where $q$ satisfies (\ref{equ:DQplus_not_unique_q_particular}). To further ensure its positive semi-definiteness, the scalars $\alpha_{1}, \ldots, \alpha_{r}$ must be the solution of the LMI problem shown in (\ref{equ:DQplus_not_unique_q_LMI}).

Similarly, if $Q$ is further specified, the matrix $F$ satisfies (\ref{equ:GivenR_Qstar_further_F}) and the vector $\halfvectorize(F - F_{0})$ must be of the form $\halfvectorize(F - F_{0}) = f + \sum_{i=1}^{r}\beta_{i}{\it \Gamma}_{i}$, where $f$ satisfies (\ref{equ:DQplus_not_unique_f_particular}) and $\beta_{1}, \ldots, \beta_{r}$ are the solution of (\ref{equ:DQplus_not_unique_f_LMI}) to guarantee the positive semi-definiteness. Thus, the solution space of $F$ is derived as (\ref{equ:GivenRQ_DFplus_not_unique}). Both of (\ref{equ:DQplus_not_unique_q_LMI}) and (\ref{equ:DQplus_not_unique_f_LMI}) are standard LMI problems, which can be solved easily with interior point method. \rulex

After analysing Problem \ref{problem:2a}, we proceed to discuss Problem \ref{problem:2b}, where only $R$ is given. Without the additional information, the problem is generally ill-posed, and the existence has been established in \citep{li2020continuous}. In Proposition \ref{proposition:P2b_onlyR_given}, we provide a sufficient condition for well-posedness of Problem \ref{problem:2b} and the analytical expressions of $Q$ and $F$ are also given in this case.  

\begin{proposition}\label{proposition:P2b_onlyR_given}
Denote the matrices $M_{Q}$ and $U_{Q}$ as 
\begin{equation}\label{equ:MQ_UQ_expression}
\begin{split}
M_{Q} & = \tilde{H}\tilde{A}^{\rm T}D_{n}(\tilde{H}_{d_{n}}^{\rm T}\tilde{H}_{d_{n}})^{-1}\tilde{H}_{d_{n}}^{\rm T}\tilde{N}_{d_{n}} + \tilde{H}_{d_{n}},    \\
U_{Q} & = \tilde{H}\tilde{A}^{\rm T}D_{n}(\tilde{H}_{d_{n}}^{\rm T}\tilde{H}_{d_{n}})^{-1}\tilde{H}_{d_{n}}^{\rm T}.
\end{split}
\end{equation}
Given $K(t)$ and $R \in \mathcal{R}_{s}$, Problem \ref{problem:2b} is well-posed if $(A, B)$ is controllable and $M_{Q}$ have full column rank, the analytical expression of  $Q$ can be uniquely recovered as:
\begin{equation}\label{equ:GivenR_only_Qstar_unique}
Q = \matrixize(D_{n}(M_{Q}^{\rm T}M_{Q})^{-1}M_{Q}^{\rm T}S_{Q}),
\end{equation}
where $S_{Q} = -U_{Q}g(t_{f}) - \tilde{H}\vectorize(G(t_{f})) + \dot{g}(t_{f})$. Then $F$ is uniquely expressed as $F = F_{0} + \tilde{F}$, where $\tilde{F}$ can be computed by (\ref{equ:Fstar_RQgiven_ABcontrol}).
\end{proposition}

\proof Based on the analysis in Theorem \ref{theorem:GivenRstar_QFstar_unique}, it follows that (\ref{equ:nsquare_derivatives_Yt_dn}) is still valid in this case. Since $\tilde{{\it \Gamma}}_{c}^{\rm T}D_{n}$ has full column rank, then $\halfvectorize(Y(t))$, where $t \in [t_{0}, t_{f}]$, can be expressed as below:
\begin{equation*}
\halfvectorize(Y(t)) = -(\tilde{H}_{d_{n}}^{\rm T}\tilde{H}_{d_{n}})^{-1}\tilde{H}_{d_{n}}^{\rm T}(\tilde{N}_{d_{n}}\halfvectorize(Q) + g(t)).
\end{equation*}
Taking the derivative of (\ref{equ:nsquare_derivatives_Yt_dn}), one can obtain that $\dot{g}(t) +\tilde{H}_{d_{n}}\halfvectorize(\dot{Y}(t)) = 0$. Then by substituting $\vectorize(\dot{Y}(t))$ and $\halfvectorize(Y(t))$, it can be rewritten as below:
\begin{equation*}
\begin{split}
\dot{g}(t) = & -\tilde{H}_{d_{n}}\halfvectorize(\dot{Y}(t)) = -\tilde{H}\vectorize(\dot{Y}(t)) \\ = & -\tilde{H}\tilde{A}^{\rm T}D_{n}\halfvectorize(Y(t)) + \tilde{H}_{d_{n}}\halfvectorize(Q) + \tilde{H}\vectorize(G(t)) \\ = & M_{Q}\halfvectorize(Q) + U_{Q}g(t) + \tilde{H}\vectorize(G(t)),
\end{split}
\end{equation*}
where $M_{Q}$ and $U_{Q}$ are shown in (\ref{equ:MQ_UQ_expression}). Then it holds that:
\begin{equation}\label{equ:MQ_SQ_unique_Q}
M_{Q}\halfvectorize(Q) = - U_{Q}g(t) - \tilde{H}\vectorize(G(t)) + \dot{g}(t).
\end{equation}
Since $M_{Q}\halfvectorize(Q)$ is constant, the right-hand side must be constant over $t \in [t_{0}, t_{f}]$. By evaluating (\ref{equ:MQ_SQ_unique_Q}) at $t = t_{f}$, it holds that $M_{Q}\halfvectorize(Q) = S_{Q}$. Since $M_{Q}$ has full column rank, then $\halfvectorize(Q)$ can be uniquely expressed as
\begin{equation*}
\halfvectorize(Q) = (M_{Q}^{\rm T}M_{Q})^{-1}M_{Q}^{\rm T}S_{Q},
\end{equation*}
which leads us to the analytical expression of $Q$ in (\ref{equ:GivenR_only_Qstar_unique}). Furthermore, with $Q$ computed as above, and given that $(A, B)$ is controllable, the matrix $F$ can also be uniquely recovered by Theorem \ref{theorem:GivenRstar_QFstar_unique}. Hence, Problem \ref{problem:2b} is well-posed, which completes the proof.   \rulex

\begin{remark}
When $M_{Q}$ does not have full column rank, let $\xi_{1}, \ldots, \xi_{s}$ span $\nullspace(M_{Q})$. The solution space of $Q$ can be characterized by $Q = \matrixize(D_{n}q + \sum_{i=1}^{s}\theta_{i}D_{n}\xi_{i})$, where $q$ satisfies $M_{Q}q = S_{Q}$ and $\theta_{1}, \ldots, \theta_{s}$ solve the standard LMI problem $\matrixize(D_{n}q) + \sum_{i=1}^{s}\theta_{i} \matrixize(D_{n}\xi_{i}) \in \mathbb{S}^{n}_{+}$. Since $\tilde{{\it \Gamma}}_{c}^{\rm T}D_{n}$ has full rank, the corresponding $F$ is unique and can be determined as shown in Theorem \ref{theorem:GivenRstar_QFstar_unique}.
\end{remark}

\section{Structural Equivalence of Two Problems}\label{sec:equivalence_p1_p2}





In this section, we further investigate the structural connection between Problem \ref{problem:one_reconstruct_R} and Problem \ref{problem:two_R_given}. By introducing an auxiliary linear time-varying quadratic regulator, we show that two IOC problems can be transformed into each other through a mapping between LQR formulations, under the assumption that $B$ and $K(t)$ have full rank for $t \in [t_{0}, t_{f}]$. The proposed methods in this paper can also be applied to inverse problems for certain class of linear time-varying quadratic regulators.

We first formulate a linear time-varying quadratic regulator as an auxiliary problem, as outlined below:
\begin{equation}\label{equ:dual_forward_problem}
\begin{split}
\min_{v} \quad & z^{\rm T}(t_{f})\bar{F}z(t_{f}) + \int_{t_{0}}^{t_{f}} \! z^{\rm T}\bar{Q}(t)z+v^{\rm T}\bar{R}(t)v \, \mathrm{d}t   \\
\text{s.t.} \quad & \dot{z}(t) = K(t)v(t), \quad z(t_{0}) = z_{0}, 
\end{split}
\end{equation}
where $z(t) \in \mathbb{R}^{m}$, $v(t) \in \mathbb{R}^{n}$. The cost matrices are given as $\bar{F} = R$, $\bar{Q}(t) = B^{\rm T}P(t)B$ and $\bar{R}(t) = P(t)$. Since $K(t)$ and $B$ have full rank for $t \in [t_{0}, t_{f}]$, it follows from Lemma \ref{lemma:necessary_sufficient_conditions_1973} that $\bar{R}(t) \succ 0$ for all $t \in [t_{0}, t_{f}]$. The matrix $\bar{Q}(t) \succ 0$ since $x_{1}^{\rm T}\bar{Q}(t)x_{1} = x_{1}^{\rm T}B^{\rm T}P(t)Bx_{1} > 0$ for all non-zero $x_{1} \in \mathbb{R}^{m}$ and $t \in [t_{0}, t_{f}]$. Then the forward optimal control problem  (\ref{equ:dual_forward_problem}) exists a unique minimizer. 

According to $B^{\rm T}P(t) = -RK(t)$ for $t \in [t_{0}, t_{f}]$, it holds that $B^{\rm T}P(t)B = -RK(t)B$ and $B = -P(t)^{-1}K(t)^{\rm T}R$. Denote $\bar{P}(t)$ as the solution of the corresponding (DRE)  of (\ref{equ:dual_forward_problem}), which can be computed as below:
\begin{equation*}
\begin{split}
\dot{\bar{P}}(t) & = \bar{P}(t)K(t)P(t)^{-1}K(t)^{\rm T}\bar{P}(t) - \bar{Q}(t) \\ & = \bar{P}(t)K(t)P(t)^{-1}K(t)^{\rm T}\bar{P}(t) + RK(t)B \\ & = \bar{P}(t)K(t)P(t)^{-1}K(t)^{\rm T}\bar{P}(t) - RK(t)P(t)^{-1}K(t)^{\rm T}R,
\end{split}
\end{equation*}
with $\bar{P}(t_{f}) = \bar{F} = R$. It follows directly that $\bar{P}(t) = R$ and the unique optimal solution for (\ref{equ:dual_forward_problem}) is given as 
\begin{equation*}
v(t) = -\bar{R}(t)^{-1}K(t)^{\rm T}\bar{P}(t)z(t) = -P(t)^{-1}K(t)^{\rm T}Rz(t) = Bz(t).  
\end{equation*}

Then consider IOC problems for (\ref{equ:dual_forward_problem}). Given the optimal feedback matrix $B$, Problem \ref{problem:one_reconstruct_R} and Problem \ref{problem:two_R_given} can also be formulated. In the next proposition, the structural equivalence of two IOC problems is characterized via a mapping between LQR formulations. For clarity, we refer to the classic and auxiliary linear quadratic regulators (\ref{equ:classic_optimal_problem}) and (\ref{equ:dual_forward_problem}) as LQR-A and LQR-B, respectively.


\begin{proposition}\label{proposition:solution_equivalence_two_problems}
Given the underlying linear systems and the corresponding optimal feedback matrices of LQR-A and LQR-B. Consider Problem \ref{problem:1a} for LQR-A with $Q$ and $F$ given, its solution is equivalent to the solution $\bar{P}(t)$ of Problem \ref{problem:2a} for LQR-B with $\bar{R}$ given. 

Similarly, consider Problem \ref{problem:2a} for LQR-A with $R$ given, the solution $P(t)$ is equivalent to the solution of Problem \ref{problem:1a} for LQR-B with $\bar{Q}$ and $\bar{F}$ given.
\end{proposition}

\proof Consider Problem \ref{problem:1a} for LQR-A with given $Q$ and $F$, the positive definite $P(t)$ of the corresponding (DRE) can be computed for all $t \in [t_{0}, t_{f}]$. Then according to Proposition \ref{proposition:GivenFQ_wellposedness}, the matrix $R$ can be uniquely recovered as (\ref{equ:unique_Rmatrix_Kt_fullrank}) for any $t_{1} \in [t_{0}, t_{f}]$ and it is the true value $R^{\star}$. Consider Problem \ref{problem:2a} on LQR-B, where the matrix $\bar{R}(t)$ is given. Then the matrix $\bar{Q}(t)$ can also be computed. By Lemma \ref{lemma:Pt_representation_by_R}, the solution $\bar{P}(t)$ of the corresponding (DRE) can be expressed as below:
\begin{equation*}
\begin{split}
\bar{P}(t) = -B^{\rm T}P(t)(P(t)BK(t))^{\dagger}P(t)B  + \bar{Y}(t),
\end{split}
\end{equation*}
where $\bar{Y}(t)$ is the zero matrix since  $K(t)^{\rm T}\bar{Y}(t) = 0$ and $K(t)$ has full row rank. By multiplying $K(t)^{\rm T}$ and $K(t)$ on left and right sides of $\bar{P}(t)$, we obtain that $K(t)^{\rm T}\bar{P}(t)K(t) = K(t)^{\rm T}R^{\star}K(t)$. Then $\bar{P}(t) = R^{\star}$ for all $t \in [t_{0}, t_{f}]$. Therefore, the solution $R$ to Problem \ref{problem:1a} for LQR-A is equivalent to the solution $\bar{P}(t)$ to Problem \ref{problem:2a} for LQR-B.


Similarly, consider Problem \ref{problem:2a} for LQR-A, when $R$ is given, the solution space of $P(t)$ is provided in \citep{jameson1973inverse}. Consider Problem \ref{problem:1a} for LQR-B with $\bar{Q}$ and $\bar{F}$ given. The positive definite solution $\bar{P}(t)$ of the corresponding (DRE) can be computed as $\bar{P}(t) = R$. By Theorem \ref{theorem:givenP_findR}, the matrix $\bar{R}(t)$ for $t \in [t_{0}, t_{f}]$ can be expressed as below
\begin{equation*}
\begin{split}
\bar{R}(t) = -K(t)^{\rm T}R(RK(t)B)^{\dagger}RK(t) +W(t)W(t)^{\rm T},
\end{split}
\end{equation*}
where $W(t)$ has $n-r_{b}$ linearly independent columns such that $B^{\rm T}W(t) = 0$. The solution space of $\bar{R}(t)$ is equivalent to the solution space of  $P(t)$ computed from Problem \ref{problem:2a} for LQR-A.  \rulex

By Proposition \ref{proposition:solution_equivalence_two_problems}, two IOC problems can be transformed into each other through a mapping between LQR formulations. It implies that inverse problems for time-varying LQR in the form of (\ref{equ:dual_forward_problem}) can be reformulated as inverse problems for time-invariant LQR, to which the methods developed in this paper can be applied.




\section{Conclusions and Future Work}\label{sec:conclusions_futurework}

In this paper, we investigated inverse optimal control problems for continuous-time LQR over finite-time horizons. Firstly, we propose two methods to recover $R$ when $Q$ and $F$ are known. The first method leverages the full trajectory of $K(t)$, for which we provide the necessary and sufficient condition for the unique recovery. When this condition holds, we derive its unique solution; otherwise, we characterize the solution space. As a complement, we introduce a second method that uses some time points of $K(t)$ to recover $R$. It can be naturally extended to the cases where only $F$ is available, and the first method fails to apply. Secondly, we studied the recovery of $Q$ and $F$ with a known $R$ in a more general manner. Unlike previous works that assume controllability a priori, our analysis clarifies its role in the well-posedness of the problem. When $F$ is specified, the solution space of $Q$ is characterized for both controllable and uncontrollable cases, and vice versa. In scenarios where both $Q$ and $F$ are unknown, we provide a sufficient condition for well-posedness and derive the unique solution. Finally, we analyse the structural equivalence between the IOC problems of reconstructing $R$ and the one in which $R$ is given as a prior. Future work will explore extensions to scenarios with partial and noisy data.



\appendix

\begin{appendices}
\section{Proof of Proposition \ref{proposition:properties_BK}}\label{appendix:proof_BK}

\proof Firstly, since $\lambda_{\tilde{t}i} < 0$, the vectors $v_{\tilde{t}i}$ are non-zero, then $Bv_{\tilde{t}i}$ are also non-zero, where $i \in [r_{\tilde{t}}]$. Then by $BK(\tilde{t}) Bv_{\tilde{t}i} = \lambda_{\tilde{t}i} Bv_{\tilde{t}i}$, we can obtain that $\lambda_{\tilde{t}i}$ are the eigenvalues of $BK(\tilde{t})$ and the corresponding eigenvectors are $Bv_{\tilde{t}i}$, where $i \in [r_{\tilde{t}}]$.

Secondly, since $\rank(K(\tilde{t})) = \rank(K(\tilde{t})B) = r_{\tilde{t}}$, the null space of $K(\tilde{t})$ contains $n - r_{\tilde{t}}$ linearly independent vectors, denoted by $w_{\tilde{t}j}$, where $j \in [n-r_{\tilde{t}}]$. Then it holds that $BK(\tilde{t})w_{\tilde{t}j} = 0$, which shows that $0$ is an eigenvalue of $BK(\tilde{t})$ and $w_{\tilde{t}j}$ are the eigenvectors, where $j \in [n-r_{\tilde{t}}]$.

Thirdly, we show that the vectors $\{Bv_{\tilde{t}i}, w_{\tilde{t}j}\}$, where $i \in [r_{\tilde{t}}]$ and $j \in [n-r_{\tilde{t}}]$, are linearly independent. We assume that there exists scalars $\alpha_{1}, \ldots, \alpha_{r_{\tilde{t}}}$, $\beta_{1}, \ldots, \beta_{n-r_{\tilde{t}}}$, which are not all zero, such that $\sum_{i = 1}^{r_{\tilde{t}}} \alpha_{i}Bv_{\tilde{t}i} + \sum_{j = 1}^{n-r_{\tilde{t}}} \beta_{j}w_{\tilde{t}j} = 0$. If $\alpha_{1}, \ldots, \alpha_{r_{\tilde{t}}}$ are all zero, then $\beta_{1}, \ldots, \beta_{n-r_{\tilde{t}}}$ are also all zero since $\{w_{\tilde{t}j}\}$ are linearly independent, where $j \in [n-r_{\tilde{t}}]$. Then both of $\alpha_{1}, \ldots, \alpha_{r_{\tilde{t}}}$ and $\beta_{1}, \ldots, \beta_{n-r_{\tilde{t}}}$ must contain non-zero scalars. Multiplying $K(\tilde{t})$, it holds that $K(\tilde{t})B(\alpha_{1}v_{\tilde{t}1} + \ldots + \alpha_{r_{\tilde{t}}}v_{\tilde{t}r_{\tilde{t}}}) = 0$. If $r_{\tilde{t}} = m$, zero is not the eigenvalue of $K(\tilde{t})B$, then the vector $\alpha_{1}v_{\tilde{t}1} + \ldots + \alpha_{r_{\tilde{t}}}v_{\tilde{t}r_{\tilde{t}}}$ must be zero. It indicates that $\alpha_{1}, \ldots, \alpha_{r_{\tilde{t}}}$ are all zero, which is the contradiction. If $1 \leq r_{\tilde{t}} < m$, we denote the eigenvectors of $K(\tilde{t})B$ corresponding with zero eigenvalue as $v_{\tilde{t}k}$, where $k-r_{\tilde{t}} \in [m-r_{\tilde{t}}]$, then $(\alpha_{1}v_{\tilde{t}1} + \ldots + \alpha_{r_{\tilde{t}}}v_{\tilde{t}r_{\tilde{t}}}) = (\sum_{k = r_{\tilde{t}}+1}^{m}\alpha_{k}v_{\tilde{t}k})$,
which contradicts with the assumption that $K(\tilde{t})B$ has $m$ linearly independent eigenvectors. Hence, the scalars $\alpha_{i}$ and $\beta_{j}$ are zero, then $\{Bv_{\tilde{t}i}, w_{\tilde{t}j}\}$, where $i \in [r_{\tilde{t}}]$ and $j \in [n-r_{\tilde{t}}]$, are linearly independent. Thus, we prove that $BK(\tilde{t})$ has $n$ linearly independent eigenvectors, which means it is diagonalizable. Then $\rank(BK(\tilde{t})) = r_{\tilde{t}}$.  \rulex


\section{Proof of Proposition \ref{proposition:solutionR_Ppositivesemi-definite}} \label{appendix:R_P_posisemi}

\proof We omit $'(\tilde{t})'$ for brevity in this proof. Denote $R^{\star}$ as the true value of the control cost matrix. It satisfies $B^{\rm T}P = -R^{\star}K$. Since $PBK$ is symmetric, it is clear that $R_{2}$ is also symmetric and it holds that:
\begin{equation*}
\begin{split}
R_{2}K = & -B^{\rm T}PK^{\dagger}K - (K^{\dagger})^{\rm T}PBK + (K^{\dagger})^{\rm T}K^{\rm T}B^{\rm T}PK^{\dagger}K \\ = & -B^{\rm T}PK^{\dagger}K - (K^{\dagger})^{\rm T}PBK + (K^{\dagger})^{\rm T}PBKK^{\dagger}K \\ = & -B^{\rm T}PK^{\dagger}K = R^{\star}KK^{\dagger}K = R^{\star}K = -B^{\rm T}P.
\end{split}
\end{equation*}
Then all the symmetric solutions are given as $R = R_{2} + Y$, where $Y = Y^{\rm T}$ and $K^{\rm T}Y = 0$. Then the set $\mathcal{V}_{R}$ in (\ref{equ:solution_set_R_positivesemi-definite_P}) follows to ensure the positive definiteness of $R$.  \rulex

\end{appendices}

\end{document}

%% file: jsscN.tex
\makeatletter
\def\@oddfoot{\hfill}
\newcount\shumeicount
\def\setshumei#1#2#3{%
  \shumeicount=\count0
  \def\@oddhead{%
    \raise-5pt\hbox to0pt{\vrule width\hsize height 0pt depth 0.4pt\hss}\relax
    \ifnum \shumeicount=\count0
      \raise-7pt\hbox to0pt{\vrule width\hsize height 0pt depth 0.4pt\hss}\relax
      #1
    \else
      \ifodd\count0
        #2
      \else
        #3
       \fi
     \fi
  }%
}
\makeatother
\makeatletter
\def\@oddfoot{\hfill}
\newcount\shujiaocount
\def\setshujiao{%
  \shujiaocount=\count0
  \def\@oddfoot{%
      \ifodd\count0
      \else
      \fi
  }%
}
\makeatother
\def\title#1#2#3#4{{
  \vspace*{0.3cm}
  \begin{flushleft} \Large\bf #1\end{flushleft}
  \vspace*{-0.2cm}
      \begin{flushleft}
      \bf #2
      \end{flushleft}
      \footnotetext{\hspace{-6mm} #3\\ #4}}}

\def\dshm#1#2#3#4
{\setshumei{J Syst Sci Complex (#1) #2:
{\thepage--\pageref{LastPage}}\hfill}
            {\hfill {\small #3}\hfill\hbox to0pt{\hss\thepage}}
            {\hbox to0pt{\thepage\hss}\hfill {\small #4}\hfill
            }
            \setshujiao}
\def\drd#1#2
{{\vskip 1cm\small \begin{flushleft}
 #1 \\
 #2\\
\copyright The Editorial Office of JSSC \&  Springer-Verlag GmbH
Germany 2021
\end{flushleft}}}

\def\hat{\widehat}
\def\tilde{\widetilde}
\def\bar{\overline}
\def\epsilon{\varepsilon}
\def\proof{\vspace{1mm}\indent {\it Proof}\quad}
